\documentclass[11pt]{article}

\usepackage{amssymb,amsmath}
\usepackage{amsthm}
\usepackage{graphicx}
\usepackage{pdfsync}
\usepackage{epstopdf}
\usepackage{lscape}
\usepackage{stmaryrd}

\usepackage{color}

\newcommand{\eol}[1]{\vspace{#1in}}
 \newtheorem{thm}{Theorem}[section]
 \newtheorem{defn}[thm]{Definition}
 \newtheorem{ex}[thm]{Example}
 \newtheorem{cor}[thm]{Corollary}
 
 \newtheorem{prop}[thm]{Proposition}
\newtheorem{rem}[thm]{Remark}
\newtheorem{obs}[thm]{Observation}

\def\G{{\mathcal G}}

\def\field{{\mathbb F}}

\def\rmnmats #1 #2{M_{#1 #2}\,({\mathbb R})}
\def\cmnmats #1 #2{M_{#1 #2}\,({\mathbb C})}

\def\mnmats #1 #2{M_{#1 #2}}

\def\matdim #1 #2{#1 \times #2}

\def\rk{{\rm rank}}

\newcommand{\sS}{\mathcal{S}^-}
\newcommand{\smr}{\operatorname{mr}^-}

\newcommand{\sMR}{\operatorname{MR}^-} 

\newcommand{\mr}{\operatorname{mr}}
\newcommand{\sM}{\operatorname{M}^-}

\newcommand{\sZ}{\operatorname{Z}^-}

\newcommand{\M}{\operatorname{M}}
\newcommand{\match}{\operatorname{match}}

\newcommand{\nty}{\operatorname{nullity}}

\newcommand{\x}{\times}
\newcommand{\bit}{\begin{itemize}}
\newcommand{\eit}{\end{itemize}}
\newcommand{\beq}{\begin{equation}}
\newcommand{\eeq}{\end{equation}}
\newcommand{\bea}{\begin{eqnarray*}}
\newcommand{\eea}{\end{eqnarray*}}
\newcommand{\bpf}{\begin{proof}}
\newcommand{\epf}{\end{proof}}

\title{Some results on minimum skew zero forcing sets, and skew zero forcing number}
\author{Luz~M.~DeAlba\thanks{Distinguished Professor Emerita, Department of Mathematics and Computer Science, Drake University, Des Moines, IA 50311, USA
(luz.dealba@drake.edu).}}

\begin{document}

\maketitle

\bigskip

\begin{abstract}

Let $G$ be a graph, and $Z$ a subset of its vertices, which we color black, while the remaining are colored white. We define the skew color change rule as follows: if $u$ is a vertex of $G$, and exactly one of its neighbors $v$, is white, then change the color of $v$ to black. A set  $Z$ is a skew zero forcing set for $G$ if  the application of the skew color change rule (as many times as necessary) will result in all the vertices in $G$ colored black. A set  $Z$ is a minimum skew zero forcing set for $G$ if  it is a skew zero forcing set for $G$ of least cardinality. The skew zero forcing number $\sZ (G)$ is  the  minimum of $|Z|$ over all  skew zero forcing sets $Z$ for  $G$. 

In this paper we discuss graphs that have extreme skew zero forcing number. We characterize complete multipartite graphs in terms of $\sZ (G)$. We note relations between minimum skew zero forcing sets and matchings in some bipartite graphs, and in unicyclic graphs.  We establish that the elements in the set of minimum skew zero forcing sets in certain bipartite graphs are the bases of a matroid.

\end{abstract}



\section{Introduction}

A {\em graph} is a pair  $G = (V_G, E_G)$, where $V_G$ is the (finite, nonempty) set of vertices of $G$ and $E_G$ is the set of edges, where an edge is a two-element subset of vertices. The {\em complete graph} on $n$ vertices is denoted $K_n$. 
An {\em induced subgraph} of $G$ is a subgraph obtained from $G$ by deleting a vertex $v$, or a number of vertices $S$, and we write $G-v$ or $G-S$, respectively. If $\{u, v \} \in E_G$ the vertices $u$ and $v$ are said to be {\em adjacent}, they are also said to be {\em neighbors}. The set $N(v)$, consisting of all the neighbors of $v$, is called the {\em open neighborhood of $v$}  (it does not include $v$), the set $N[v] = N(v) \cup \{v \}$ is the {\em closed neighborhood of $v$}. The {\em degree of a vertex $v \in V_G$}, denoted by $\deg_G (v)$, is the number of edges adjacent to $v$. The minimum (respectively, maximum) degree in a graph $G$ is denoted $\delta (G)$ (respectively, $\Delta (G)$).  A subset $S \subseteq V_G$ is called {\em independent} if no two vertices in $S$ are adjacent. 
A graph $G$ is {\em connected} if each pair of vertices in $V_G$ belongs to a path. A vertex $v \in V_G$ is a {\em cut-vertex} if the induced graph $G-v$ is not connected. We say that $G$ is the  {\em vertex sum} of two graphs $G_1$ and $G_2$, and write $G_1 \bigoplus_v G_2$ if $v$ is a cut-vertex of $G$, $V_{G_1} \cap V_{G_2} = \{v \}$, and $E_{G_1} \cap E_{G_2} = \emptyset$. A graph with no cut-vertices is said to be {\em nonseparable}. 

A {\em matching} in a graph $G$ is a set of edges $M = \{ \{ i_1, j_1\}, \{ i_2, j_2 \}, \dots, \{ i_k, j_k \} \} \subseteq E_G$, such that no endpoints are shared. The vertices that determine the edges in $M$ are called {\em $M$-saturated} vertices, all other vertices in $V_G$ are called {\em $M$-unsaturated} vertices. A {\em perfect matching} in a graph $G$ is a matching that saturates all vertices of $G$. A {\em maximum matching} in a graph $G$ is a matching of maximum order among all matchings in $G$. The {\em matching number} of a graph $G$, denoted by $\match(G)$, is the number of edges in a maximum matching. An even  cycle in a graph $G$ is called {\em $M$-alternating} if it alternates between edges in $M$ and edges not in $M$. A matching $M$ in a graph $G$ is {\em uniquely restricted} if $G$ does not contain an $M$-alternating cycle.

A graph $G$ is {\em $k$-partite} if $V_G$ can be expressed as the union of $k$ (possibly empty) independent sets, and is denoted $K_{n_1,n_2,\dots,n_k}, k \ge 2, n_i \ge 1, i = 1, 2, \dots, k$. A {\em tree} is a connected graph $T$, with $\left| E_T \right| = \left| V_T-1 \right|$, trees are {\em 2-partite}, also known as {\em bipartite}.

Although many of the results presented here are valid for some finite fields, we assume throughout this paper that $\field$ is an infinite field. A matrix $A \in \field^{n \x n}$ is {\em skew-symmetric} if $A^T = -A$. For an $n \x n$ skew-symmetric matrix $A$, the {\em graph of $A$}, denoted $\G(A)$, is the graph with vertices $\{v_1, . . . , v_n \}$ and edges $\{ \{v_i, v_j \} : a_{ij} \ne 0, 1 \le i < j \le n \}$.

Let $\sS(\field, G) = \{ A \in \field^{n \x n} : A^T = -A, \G(A) = G \}$ be the set of  skew-symmetric matrices over the field $\field$ described by a graph  $G$. The {\em minimum skew rank of a graph $G$ over the field $\field$ } is defined as $\smr(\field, G) = \min \{ \rk(A) : A \in \sS(\field, G) \}$, the {\em maximum skew nullity of $G$ over the field $\field$} is defined as $\sM(\field, G) = \max \{\nty (A) : A \in \sS(\field, G) \}$, and the {\em maximum skew rank of $G$ over the field $\field$} as $\sMR(\field, G) = \max \{ \rk(A) : A \in \sS(\field, G) \}$. Clearly $\smr(\field, G) + \sM(\field, G) = |G|$, but note that, since a skew symmetric matrix has even rank, $\sMR (\field, G) \le |G|$.

For a graph $G$,  select $Z \subseteq V_G$, color all vertices in $Z$ black, and all others white. Next apply the {\em skew color change rule}: if $u \in V_G$ ($u$ any color), and exactly one of its neighbors $v$, is white, then change the color of $v$ to black (we say $u$ forces $v$ black). Continue to  apply the skew color change rule until no more changes are possible. A {\em skew zero forcing set} for a graph $G$ is a subset  $Z$ of $V_G$, such that, if initially the vertices in $Z$ are colored black and the remaining vertices are colored white, the skew color change rule forces all the vertices in $V_G$ black. A {\em minimum skew zero forcing set} for a graph $G$ is a skew zero forcing set of minimum order among all skew zero forcing sets for $G$. The {\em skew zero forcing number $\sZ(G)$} is the minimum of $|Z|$ over all skew zero forcing sets $Z \subseteq V_G$.

\section{Preliminary results}\label{prelim}

The parameter $Z(G)$ was introduced in~\cite{AIM}, while the parameter  $\sZ(G)$, was introduced in~\cite{10IMA}. 

\begin{prop}\label{tree-known}
\begin{enumerate}
\item\label{induced}\cite[Observation~1.7]{10IMA}
If $H$ is an induced subgraph of $G$, then $\smr (\field, H) \le \smr (\field, G)$.
\item\label{2.810IMA}~\cite[Proposition 3.5]{10IMA} For any graph $G$, $\sM (\field, G) \le \sZ (G)$ and $\smr (\field, G) \ge |G| - \sZ (G)$.
\item\label{M=Z}~\cite[Proposition 4.2]{AIM} For any tree $T$, $M (\field, T) = Z (T)$, and hence $\mr (\field, T) = |T| - \M (\field, T) = |T| - Z (T)$.
\end{enumerate}
\end{prop}

\begin{thm}\label{knowngen}
\begin{enumerate}
\item\label{smr=2}\cite[Theorem~2.1]{10IMA} Let $G$ be a connected graph with $|G| \ge 2$, then $\smr (\field, G) = 2$ if and only if $G = K_{n_1,n_2,\dots,n_s}, s \ge 2, n_i \ge 1, i = 1, 2, \dots, s$.
\item\label{mr=match}~\cite[Theorem~2.5]{10IMA}  For a graph $G$, $\sMR (\field, G) = 2 \match (G)$, and every even rank between $\smr (\field, G)$ and $\sMR (\field, G)$ is realized by a matrix in $\sS(\field, G)$.
\item\label{upm}~\cite[Theorem~2.6]{10IMA} For a graph $G$, $\smr(\field, G) = |G| = \sMR(\field, G)$ if and only if $G$ has a unique perfect matching.
\item\label{smrtree}\cite[Theorem~2.8]{10IMA}
If $T$ is a tree, then $\smr(\field, T) = 2 \match(T) = \sMR(\field, T)$.
\item\label{cover}\cite[Proposition 3.3]{10IMA}. Let $\field$ be a field and $G = \cup_{i=1}^k \, G_i$ be a graph. Suppose that for all $i \ne j$, $G_i$ and $G_j$ have no edges in common, then $\smr (\field, G) \le \sum_{i=1}^k \, \smr ( \field, G_i)$.
\end{enumerate}
\end{thm}

\section{Graphs with extreme skew zero forcing number}\label{extreme}

It is a fact that for any graph $G$, $0 \le \sZ (G) \le |G|$. If a graph has isolated vertices, those vertices must belong to all skew zero forcing sets for the graph. Thus, without loss of generality, we assume that graphs have no isolated vertices. Also, some of the results presented here are valid for graphs that are disconnected,  we specifically note when a graph must be connected.

\begin{rem}\label{zsub} If $G$ is a graph, $v \in V_G$, and $Z$  a minimum skew zero forcing set for $G-v$, then $Z \cup \{v \}$ is a skew zero forcing set for $G$, so $\sZ(G) \le \sZ(G - v) +1$.
\end{rem}

From Proposition~\ref{tree-known}, and Theorem~\ref{knowngen}, we obtain the  inequalities 

\begin{equation}\label{ineq}
|G| - \sZ (G) \le \smr (\field, G) \le \sMR (\field, G) = 2 \match (G) \le |G|.
\end{equation}

The following are derived using the inequalities in Equation~\ref{ineq}, and the definition of  skew zero forcing set. 

\begin{obs}\label{smallz}
\begin{enumerate}
\item\label{z=0} If $\sZ (G) = 0$, then there is a vertex, $v \in V_G$, such that $\deg_G (v) = 1$.
\item\label{z=0upm} If $\sZ (G) = 0$, then $|G|$ is even, $\smr (\field, G) =  |G|$ and $G$ has a unique perfect matching.
\item\label{z=1even} If $\sZ (G) = 1$ and $|G|$ is even, then $\smr (\field, G) =  |G|$ and $G$ has a unique perfect matching.
\item\label{z=1odd}  If $\sZ (G) = 1$ and $|G|$ is odd, then $\smr (\field, G) =  |G| - 1.$
\item\label{z=2odd}  If $\sZ (G) = 2$ and $|G|$ is odd, then $\smr (\field, G) =  |G| - 1.$
\item\label{z=2even} If $\sZ (G) = 2$ and $|G|$ is even, then either $\smr (\field, G) =  |G|$ and $G$ has a unique perfect matching, or $\smr (\field, G) =  |G| -2$.
\item (\cite[Observation 1.6]{10IMA} $\sZ (G) = |G|$ if and only if $G$ consists only of isolated vertices.
\end{enumerate}
\end{obs}

\begin{figure}[h!]
\begin{center}
\scalebox{.4}{\includegraphics{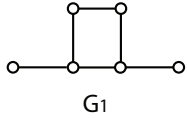}}
\caption{Graph with $\sZ (G_1) = 0 = |G_1|-6$.}\label{upmsz0}
\end{center}
\end{figure}

\begin{figure}[h!]
\begin{center}
\scalebox{.4}{\includegraphics{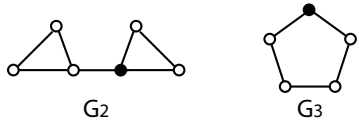}}
\caption{Graphs with $\sZ (G_i) = 1 = |G_2| - 5 = |G_3| - 4$.}\label{upmsz1}
\end{center}
\end{figure}

\begin{figure}[h!]
\begin{center}
\scalebox{.4}{\includegraphics{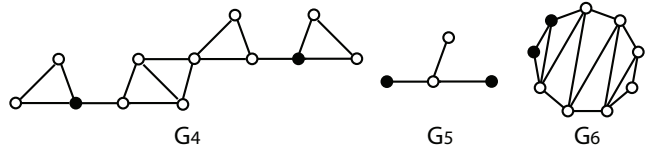}}
\caption{Graphs with $\sZ (G_i) = 2 = |G_4| -10 = |G_5| -2 = |G_6| - 7$.}\label{upmsz2}
\end{center}
\end{figure}

It is clear that the converses of Items~\ref{z=0}--\ref{z=2even} in Proposition~\ref{smallz}  are not true, the graphs in Figures~\ref{upmsz0}--\ref{upmsz2}, and Item~\ref{tri} in Observation~\ref{special}, also illustrate this. For the graph $G_2$ in Figure~\ref{upmsz1} (which is a cactus graph, and also a block-clique graph), $\smr (\field, G) =  |G| > |G| - 1 = |G| - \sZ (G)$. 

\eol{.125}

Note that if $G$ is one of the graphs $K_2$, $K_3$, or $K_{2,1}$, then $\sZ (G) = |G| - 2$. We now show that this equation characterizes all complete multipartite graphs. The  proof given below involves the use of $\smr (\field, G)$, but one can easily construct a field independent proof.

\begin{thm}\label{sz=g-2}
A connected graph $G$ is a complete multipartite graph $K_{n_1, n_2, \dots, n_s}$, $s \ge 2, n_i \ge 1$ if and only if $\sZ(G) = |G|-2$.
\end{thm}

\bpf
Let $G = K_{n_1, n_2, \dots, n_s}$, $s \ge 2, n_i \ge 1$, with $|G| \ge 4$. From Item~\ref{smr=2} in Theorem~\ref{knowngen}, $\smr (\field, G) = 2$, hence from Equation~\ref{ineq}, $|G| - 2 \le \sZ (G)$. Pick $u, v \in V_G$  adjacent (in different partite classes), then $Z = V_G - \{u, v \}$ is a skew zero forcing set for  $G$, and $\sZ(G) \le |G| - 2$. It follows that $\sZ(G) = |G| - 2$.

Conversely, if $G$ is connected, but not $K_{n_1, n_2, \dots, n_s}$, $s \ge 2, n_i \ge 1$, then $G$ has an induced $P_4$ or an induced paw (\cite[Remark 2.2]{10IMA}). If $v_1, v_2, v_3, v_4$ induce a $P_4$, or the paw, then $Z = V_G - \{v_1, v_2, v_3, v_4 \}$ is a skew zero forcing set for $G$, so $\sZ(G) \le |G| - 4 \ne |G| - 2$.
\epf

\begin{cor}\label{g-mr}
If $G$ is a connected graph, then either $\sZ (G) = |G| -2$ or $\sZ (G) \le |G| -4$. There are no connected  graphs for which $\sZ (G) = | G|-1$, and no connected  graphs for which $\sZ (G) = |G|-3$. 
\end{cor}

\begin{rem}\label{ord8}
One can verify, directly, that for connected graphs $G$ (pictured in~\cite{98RW}), with $4 \le |G| \le 6$, $\sZ(G) = |G| - 4$ if and only if $\smr (\field, G) = 4$. Henceforth we assume $|G| \ge 7$.
\end{rem}

\begin{prop}\label{cutv} 
If $G$  is connected, has a cut-vertex, and $\sZ (G) = |G| - 4$, then $\smr (\field, G) = 4$.
\end{prop}

\bpf

Let $G$ be connected, and $v \in V_G$ be a cut-vertex. Let $G_1$ be the connected subgraph of $G$ induced by the vertices of one of the components of $G-v$ and $v$, and $G_2$ be the connected subgraph of $G$ induced by $(V_G - V_{G_1}) \cup \{ v\}$, so that $G= G_1 \bigoplus_v G_2$.

If $\sZ (G_1) = |G_1| - 2$, and $\sZ (G_2) = |G_2| - 2$, then $G$ is the vertex sum of two complete multipartite  graphs, and in this case $\smr (\field, G) = 4$. 

The two other possibilities that arise from Corollary~\ref{g-mr} do not allow $\sZ (G) = |G| - 4$. Let $Z_1$ and $Z_2$ be minimum skew zero forcing sets for $G_1$ and $G_2$, respectively. 

If $\sZ (G_1) = |G_1| - 2$,  $\sZ (G_2) \le |G_2| - 4$, and $Z_1 \ne \emptyset$, then from the proof of Theorem~\ref{sz=g-2}, we can take $v \in Z_1$, thus  $Z_1 \cup Z_2$ is a skew zero forcing set for $G$. If $\sZ (G_i) \le |G_i| - 4, i = 1, 2$, then $Z_1 \cup Z_2  \cup \{ v \}$ is a skew zero forcing set for $G$. In both cases, $\sZ (G) < |G| - 4$.
\epf

\begin{prop}\label{induce6}
If $G$ is a connected graph, and $H$ is  a connected induced subgraph of $G$ of order 6 that has a unique perfect matching, then $\sZ (G ) < |G| -4$.
\end{prop}

\begin{figure}[h!]
\begin{center}
\scalebox{.5}{\includegraphics{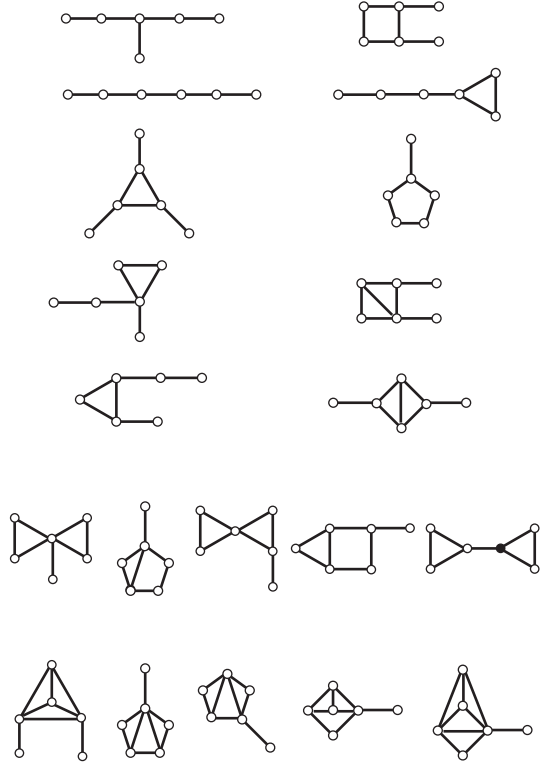}}
\caption{Graphs on six vertices~(\cite{98RW}) with a unique perfect matching.}\label{upm1}
\end{center}
\end{figure}

\bpf
Figure~\ref{upm1} shows the twenty connected graphs on six vertices that have a unique perfect matching. One of these is $G_2$, also pictured in Figure~\ref{upmsz1}, and satisfies $\sZ (G_2 ) = 1$, all others have $\sZ (H ) = 0$. If $H=G_2$, and $u \in V_{G_2}$, then $(V_G-V_H) \cup \{ u \}$ is a  skew zero forcing set of $G$; if $H \ne G_2$, then $V_G - V_H$ is a  skew zero forcing set for $G$.  Thus, $\sZ (G ) \le |G| -5 < |G| - 4$. 
\epf

\begin{prop}\label{g-4a}
Let $\field$ be a field,  and $G$ a connected graph with $|G| \ge 4$. If $\smr (\field, G) =4$, then $\sZ(G) = |G| - 4$. 
\end{prop}

\bpf
If $\smr (\field, G) = 4$, then $\sZ (G) \le |G| - 4$, so $4 \le |G| - \sZ (G) \le \smr (\field, G) =4$.
\epf

The following example provided by Sudipta Mallik and Bryan Shader, and constructed using their methods as in~\cite{13MS}, shows that the converse of Proposition~\ref{g-4a} is not true. 

\begin{defn}\cite[p. 3651]{13MS}
A collection $\{ N_i : i \in \mathcal{I} \}$ of vectors is a minimally dependent set of vector if it is a linearly dependent set and for each $j \in \mathcal{I}, \{ N_i : i \ne j, i \in \mathcal{I} \}$ is a linearly independent set of vectors.
\end{defn}

\begin{ex}
If $G =K_3 \times K_3$ is the graph with adjacency matrix and graph as in Figure~\ref{counterex}, 

\begin{figure}[h!]
\begin{center}
\scalebox{.5}{\includegraphics{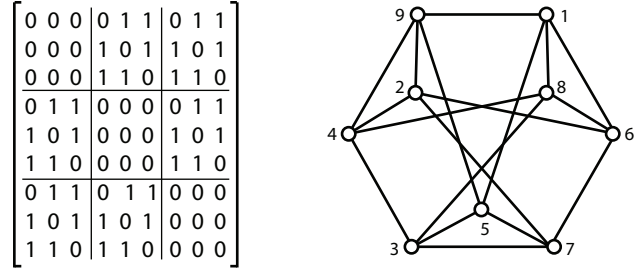}}
\caption{The graph $K_3 \times K_3$.}\label{counterex}
\end{center}
\end{figure}

then $\sZ (G) = |G| -4$, and $\smr (G) \ge 6$.
\end{ex}

\bpf
Using the graph in Figure~\ref{counterex}, and the fact that $G$ is tripartite, it is not difficult to show that $\sZ(G) > 4$, and that $Z= \{ 1,2,3,4,7 \}$ is a minimum skew zero forcing set for $G$. Hence $\sZ (G) = 5 = |G|-4$, and  $4 = |G| - \sZ(G) \le \smr (G)$. 

Let $B \in \sS (G)$, and assume columns $3i + 1, 3i + 2, 3i + 3$ are linearly independent for $i = 0, 1, 2$. Then from the zero-nonzero pattern of $B$ we observe that columns $3i + 1, 3i + 2, 3i + 3, 3i + 4, 3i + 5 \ (\mbox{mod } 9)$ are linearly independent, and since $B$ is skew symmetric,  $\rk (B) \ge 6$. 

Assume now that columns $3i + 1, 3i + 2, 3i + 3$ of $B$ are linearly dependent for $i = 0, 1, 2$, and hence minimally linearly dependent. By Lemma 4.7 in~\cite{13MS}, the nullspace of $B$ contains vectors of
the form
\begin{equation}
\left[
\begin{array}{c}
a\\
b\\
c\\
0\\
0\\
0\\
0\\
0\\
0
\end{array}
\right],
\left[
\begin{array}{c}
0\\
0\\
0\\
d\\
e\\
f\\
0\\
0\\
0
\end{array}
\right], \mbox{ and }
\left[
\begin{array}{c}
0\\
0\\
0\\
0\\
0\\
0\\
g\\
h\\
k
\end{array}
\right],\end{equation}\label{vecnul}

for some $a, b,  c,  d,  e,  f,  g,  h,  k$, each of which is nonzero. Let $D = \mbox{diag} (a, b,  c,  d,  e,  f,  g,  h,  k)$, hence $DBD \in \sS(G)$, and $\rk (DBD) = \rk (B)$. Note that the nullspace of $DBD$ contains vectors as in Equation~\ref{vecnul} with $a = b =  c =  d =  e =  f =  g =  h =  k = 1$. Direct calculations now show that $DBD$ has the form: $$\left[
\begin{array}{ccc}
0 & x S & y S\\
x S & 0 & z S\\
y S & z S & 0
\end{array}
\right] =
\left[
\begin{array}{ccc}
0 & x & y\\
x & 0 & z\\
y & z &0
\end{array}
\right] 
\otimes
S$$ for some nonzero $x, y$ and $z$. Since $\rk \left( \left[
\begin{array}{ccc}
0 & x & y\\
x & 0 & z\\
y & z &0
\end{array}
\right] \right) = 3$, and $\rk (S) = 2$,  it follows that $\rk (B) = \rk (DBD) = 6$. Hence $\smr(G) \ge 6$.
\epf

\section{Bipartite graphs}\label{main}

In this section we study the relation between certain matchings and skew zero forcing sets.  Bipartite graphs provide a natural setting for this discussion.

\begin{prop}\label{constrained}
If $B$ is a  bipartite graph, and $M$ a uniquely restricted matching in $B$, then the set of $M$-unsaturated vertices of $B$ is a  skew zero forcing set for $B$. 
\end{prop}

\bpf
Let $B$ be a bipartite graph, $M$ a uniquely restricted matching in $B$, and $H$ the connected subgraph of $B$ induced by the vertices in $M$ (If $H$ is not connected, the following process  can be applied separately to each of the components of $H$). Suppose the vertices in the bipartition of $H$ are $u_1, \dots, u_r$ and $v_1, \dots, v_r$,  $\{u_i, v_i \} \in M$, $\{u_i, v_j \} \notin E_H$ whenever $1 \le i < j \le r$. Let $Q = V_B - V_H$, and color the vertices in $Q$ black. Without loss of generality we may assume $\deg_H (v_r) = 1$, then we have the following sequence of forces $v_r \rightarrow u_r, v_{r-1} \rightarrow u_{r-1}, \dots, v_1 \rightarrow u_1, u_1 \rightarrow v_1, u_2 \rightarrow v_2, \dots, u_r \rightarrow v_r$.  Thus $Q$  forms a skew zero forcing set for $B$. 
\epf

\begin{prop}\label{qlessz}
Let $G$ be a  graph,  $M$ a matching in $G$,  and $\smr (\field, G) \le 2|M|$. If the set of $M$-unsaturated vertices of $G$ is a skew zero forcing set for $G$, then it is a minimum skew zero forcing set for $G$. 
\end{prop}

\bpf
Let $\field$ be a field, $M$ a matching in $G$, and $Q$ the set of $M$-unsaturated vertices. From Equation~\ref{ineq}, $|G| - \sZ (G) \le \smr (\field, G) \le 2 |M|$, so $|Q| = |G| -  2 |M| \le  \sZ (G)$. Thus, if $Q$ is a zero forcing set for $G$, it is a minimum skew zero forcing set for $G$. 
\epf

\begin{cor}\label{z=q}
Let $G$ be a  graph,  $M$ a matching in $G$,  and $\smr (\field, G) \le 2|M|$. If the set of $M$-unsaturated vertices is a minimum skew zero forcing set for $G$, then $|G| - \sZ (G) = \smr (\field, G)$. 
\end{cor}

\bpf
If $Q$ denotes the set of $M$-unsaturated vertices, then $2 |M| = |G| - |Q| = |G| - \sZ (G) \le \smr (\field, G) \le 2|M|$.
\epf

\begin{cor}\label{bitt}
If $B$ is a bipartite graph (in particular a tree, or cactus) with $\smr (\field, B) = \sMR (\field, B)$, then 

\begin{enumerate}
\item\label{maxmr} There is a maximum matching in $B$ such that the set of $M$-unsaturated vertices is a minimum skew zero forcing set for $B$.
\item\label{mrt=t-z} $\sZ (B) = |B| - \smr (\field, B)$, and $\sM (\field, B) =  \sZ (B)$. 
\item $\sZ (B) =0$, if and only if $B$ has a unique perfect matching.
\item $\sZ (B) \le \frac{ \Delta (B) |B| - 2 |E_B|}{\Delta (B)}$. In particular, 
\begin{enumerate}
\item if $T$ is a tree, $\sZ (T) \le \frac{|T| ( \Delta (T) -2) + 2}{\Delta (T)}$,  and this bound is sharp for paths and stars;
\item if $U$ is a unicyclic, $\sZ (U) \le \frac{|U| ( \Delta (U) -2)}{\Delta (U)}$.
\end{enumerate}
\end{enumerate}
\end{cor}

\bpf
\begin{enumerate}
\item If  $B$ is a bipartite graph with $\smr (\field, B) = \sMR (\field, B)$, then there must be a uniquely restricted maximum matching $M$,  in $B$. Then use Proposition~\ref{constrained}, and  Proposition~\ref{qlessz}.
\item\label{m=z} This follows from Item~\ref{maxmr} above, and Corollary~\ref{z=q}.
\item This follows from Item~\ref{upm} in Theorem~\ref{knowngen}, as well as Item~\ref{mrt=t-z} above.
\item This follows from the fact that a bipartite graph $B$ has a matching of size at least $\frac{|E_B|}{\Delta (B)}$ (\cite[p. 108]{08HHM}).
\end{enumerate}
\epf

The results in Proposition~\ref{constrained} and Proposition~\ref{qlessz}, suggest a certain duality between maximum matchings  and minimum skew zero forcing sets. We establish that a natural duality, via matroids,  does indeed exist for some families of graphs.

\begin{defn}{\bf The matching matroid and its dual.}~\cite[pp. 92--93]{86LP}
If $G$ is a bipartite graph, the set $$\mu =\{ X \subseteq V_G \, : \, X \, \mbox{  is saturated by some matching} \}$$ is a matroid on $V_G$ with bases the sets of vertices saturated by maximum matchings in $G$. Its dual matroid $$\mu^* =\{ Q \subseteq V_G \, : \, Q \, \mbox{  is not saturated by some maximum matching} \}$$ on $V_G$ has bases $V_G - B_i$, where $B_i$ is a basis of $\mu$.
\end{defn}

\begin{thm}\label{matroid}
{\bf The zero forcing matroid.} In a bipartite graph $B$,  in which all maximum matchings are uniquely restricted, the elements in the set of minimum skew zero forcing sets are the bases of a matroid on the vertices of the corresponding graph, and this matroid is the dual of the matching matroid on $B$.
\end{thm}

\begin{thm}\label{treeall}
If $B$ is a bipartite graph in which all maximum matchings are uniquely restricted, and $M$ a matching in $B$, then $M$ is a maximum matching in $B$ if and only if the set of $M$-unsaturated vertices of $B$ is a minimum skew zero forcing set for $B$. Alternatively, let $Z$ be a skew zero forcing set for $B$, then $Z$ is a minimum skew zero forcing set for $B$ if and only if $V_B - Z$ has a unique perfect matching which is a maximum matching in $B$.
\end{thm}

\bpf
Let $B$ be a a bipartite graph in which all maximum matchings are uniquely restricted, $M$ a maximum matching in $B$, and $Q$ the set of $M$-unsaturated vertices. Since $M$ is a uniquely restricted matching, from Proposition~\ref{constrained}, $Q$ is a minimum skew zero forcing set for $B$. Conversely, suppose that $Q$ is a minimum skew zero forcing set for $B$. From Theorem~\ref{matroid}, $V_B-Q$ has a unique perfect matching which is a maximum matching in $B$.

We omit the alternate proof.
\epf

In a tree all maximum matchings are uniquely restricted matchings, thus we have the following.

\begin{cor}\label{tree}
If $T$ is a tree, and $M$ a matching in $B$, then $M$ is a maximum matching in $T$ if and only if the set of $M$-unsaturated vertices of $T$ is a minimum skew zero forcing set for $T$. Alternatively, let $Z$ be a skew zero forcing set for $T$, then $Z$ is a minimum skew zero forcing set for $T$ if and only if $V_T - Z$ has a unique perfect matching which is a maximum matching in $T$.
\end{cor}

\section{Unicyclic Graphs}

Results on the minimum skew rank of unicyclic graphs can be found in~\cite{11D}, explicitly:  $\smr (\field, U) = \sMR(\field, U)$ if the unique cycle is odd, or if the unique cycle is even and $U$ has a uniquely restricted maximum matching; $\smr (\field, U) = \sMR(\field, U) - 2$ if the unique cycle is even and $U$ does not have a uniquely restricted maximum matching.

\begin{prop}\label{uniodd}
If $U$ is a unicyclic graph, then there exists a matching $M$,  in $U$, such that the set of $M$-unsaturated vertices of $U$ is a minimum skew zero forcing set for $U$.
\end{prop}

\bpf
If the unique cycle has odd order, and $M$ is a maximum matching in $U$, we can construct a proof by induction to show that the set of $M$-unsaturated vertices is a minimum skew zero forcing set for $U$. The base cases follow from examples in~\cite{11D}; we omit the details of the proof.

If the unique cycle has even order, and $\smr (\field, U) = \sMR(\field, U)$, then the result follows from Item~\ref{maxmr} in Corollary~\ref{bitt}. 

If the unique cycle has even order, $\smr (\field, U) = \sMR(\field, U)-2$, and $\widehat M$ is a maximum matching in $U$, then the cycle is $\widehat M$-alternating. If $U$ is a cycle, and $e $ is an edge in $\widehat M$, then $M = \widehat M-e$ is a uniquely restricted matching in $U$.

If $U$ is not a cycle, then it has an induced subgraph $H$, consisting of the vertex sum of the cycle and a path of order 3, that is $H= C \bigoplus_{v_1} P_3$, where $P_3 = (\{v_1, v_2, v_3\}, \{\{v_1, v_2 \}, \{ v_2 , v_3 \} \})$, and $u$ is a neighbor of $v_1$ on the cycle. Thus, there exists a maximum matching $\widehat M$, in $U$, containing the edges $\{ u, v_1 \}$, and $\{ v_2, v_3 \}$ (if $v_2$ is $\widehat M$-unsaturated, then $(\widehat M - \{ \{ u, v_1 \} \} ) \cup \{ \{ v_1, v_2 \}\}$ is a uniquely restricted maximum matching). The matching $M = ( \widehat M -  \{ \{ u, v_1 \}, \{ v_2, v_3 \} \} ) \cup \{ \{ v_1, v_2 \} \}$, is a uniquely restricted matching in $U$. 

In either case $M$ has order $\frac{\sMR(\field, U) - 2}{2}$, and since $\smr (U) = \sMR(\field, U) - 2$, it follows from Proposition~\ref{constrained} that  the set of $M$-unsaturated vertices is a minimum skew zero forcing set for $U$. 
\epf

\begin{cor}\label{uni}
If $U$ is a unicyclic graph, then $\sZ (U) = |U| - \smr (\field, U)$.
\end{cor}

\section{Additional Examples}

We conclude with several contrasting examples of graphs $G$, for which there is a matching $M$, in $G$, such that the set of $M$-unsaturated vertices is a minimum skew zero forcing set for $G$. Also, in Observation~\ref{special}, we list the skew zero forcing number of some special graphs.

\begin{figure}[h!]
\begin{center}
\scalebox{.4}{\includegraphics{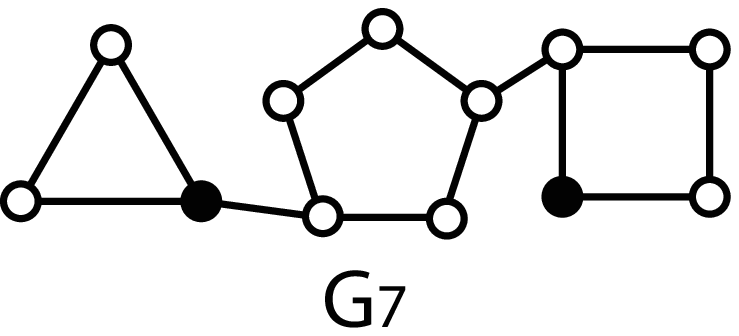}}
\caption{Non-bipartite  graph with $\sZ (G_7) = 2, |G_7| - 2 = \smr (\field, G_7)$.}\label{cactus}
\end{center}
\end{figure}

\begin{ex}
The graph $G_7$, in Figure~\ref{cactus}, is a non-bipartite cactus, does not have a unique maximum matching, but has a maximum matching, $M_7$, of order 5, with no $M_7$-alternating cycle. The set of $M_7$-unsaturated vertices (in black) is a minimum skew zero forcing set for $G_7$, and from Item~\ref{z=2even} in Observation~\ref{smallz}, $|G_7| - \sZ (G_7) = 10 = \smr (\field, G_7) \ne \sMR (\field, G_7) = 12$.  
\end{ex}

\begin{figure}[h!]
\begin{center}
\scalebox{.4}{\includegraphics{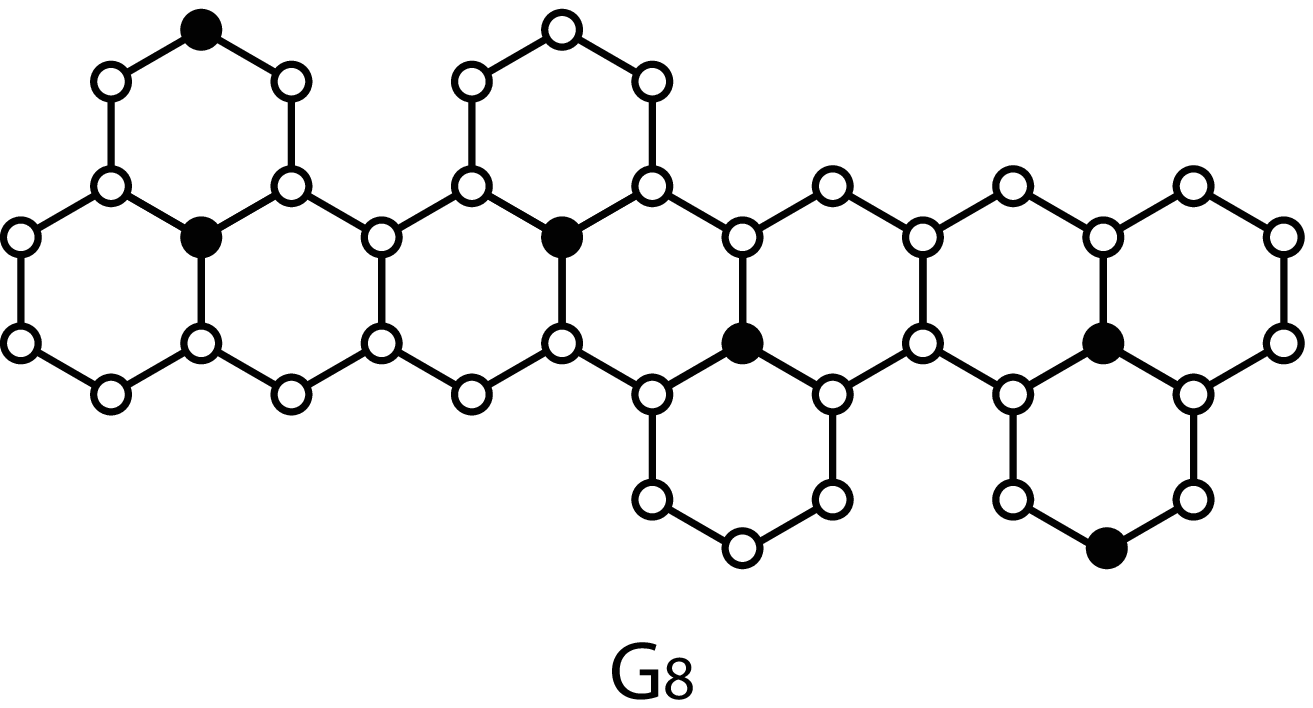}}
\caption{Bipartite  graph with $\sZ (G_8) = 6, |G_8| - 6 = \smr (\field, G_8)$.}\label{lovasz}
\end{center}
\end{figure}

\begin{ex}
The graph $G_8$ (see~\cite[pp.6--7]{86LP}), in Figure~\ref{lovasz}, is bipartite, does not have a perfect matching, but has a uniquely restricted matching $M_8$, of cardinality 18 (it is easy to verify that all matchings in $G_8$ of cardinalities 20 and 19 are not uniquely restricted), the set of  $M_8$-unsaturated vertices (in black) is a minimum skew zero forcing set for $G_8$, and with the aid of Mathematica one can verify that $|G_8| - \sZ (G_8) = 42 - 6 = 36 = \smr (\field, G) \ne \sMR (\field, G)$.
\end{ex}

\begin{ex}
The graph $G_4$, in Figure~\ref{upmsz2} is not bipartite, not unicyclic, and not a cactus, has a unique perfect matching, but also has a matching $M_4$ of order $5$, such that the set of $M_4$-unsaturated vertices is a minimum skew zero forcing set for $G_4$, and $|G_4| - \sZ (G_4) = 10 < \smr (\field, G) = \sMR (\field, G)$.
\end{ex}

Below we list some graphs and their skew zero forcing numbers. We refer the reader to~\cite{10IMA} for the definitions of $W_n$, the wheel on $n$ vertices; $P_{m,k}$, the $m, k$-pineapple, with $m \ge 3, k \ge 1$;  $Q_s$, the $s$th hypercube; $T_n$, the super-triangle; $H_s$ the $s$th half-graph; $N_s$, the necklace with $s$ diamonds; $G \circ H$, the corona of $G$ with $H$; $G \boxempty H$, the Cartesian product of $G$ and $H$.

\begin{obs}\label{special}
For the  graphs $G$ in Items~\ref{path},~\ref{cycle},~\ref{wheel},~\ref{pineapple},~\ref{cube},~\ref{half},~\ref{corona}, and~\ref{pdp}, $\sZ(G) = |G| -  \smr (\field, G)$ (note that there might be  restrictions on the field $\field$, see~\cite{10IMA}):

\begin{enumerate}
\item\label{path} $\sZ (P_n) = \left\{
\begin{array}{lll}
0 & \mbox{if} & n \mbox{ is even},\\
1 & \mbox{if} & n \mbox{ is odd};\\
\end{array}
\right.$
\item\label{cycle} $\sZ (C_n) = \left\{
\begin{array}{lll}
2 & \mbox{if} & n \mbox{ is even},\\
1 & \mbox{if} & n \mbox{ is odd};\\
\end{array}
\right.$
\item\label{wheel} $\sZ (W_n) = \left\{
\begin{array}{lll}
2 & \mbox{if} & n \mbox{ is even},\\
3 & \mbox{if} & n \mbox{ is odd};\\
\end{array}
\right.$
\item\label{pineapple} $\sZ (P_{m,k} ) = |P_{m,k} | - 4= m+k-4, m \ge 3, k \ge 1 $;
\item\label{cube} $\sZ (Q_s) = 2^{s-1}, s \ge 2$;
\item\label{tri} $\sZ (T_n) = n-1$;
\item\label{half} $\sZ (H_s ) = 0$;
\item\label{neck} $\sZ (N_s ) = s$, $\sZ(N_s) = |N_s| -  \smr (\field, N_s)$ if and only if $s=2$;
\item\label{corona} $\sZ (G \circ K_1 ) = 0$;
\item\label{corona2} $\sZ (C_t \circ K_s ) = st-3t+2, s \ge 2,  \sZ(C_t \circ K_s ) = |C_t \circ K_s | -  \smr (\field, C_t \circ K_s )$ if and only if $s$ is even;
\item\label{pdp} $\sZ (P_s \square P_s ) = s$; 
\item\label{cdp} $\sZ (K_3 \square P_2 ) = 2$, $\sZ(K_3 \square P_2  ) = |K_3 \square P_2  | -  \smr (\field, K_3 \square P_2  )$, and for $s \ge 3, t \ge 3$, $\sZ (K_s \square P_t ) = s$, $\sZ(K_s \square P_t  ) = |K_s \square P_t  | -  \smr (\field, K_s \square P_t  )$ if and only if $s$ is even, or both $s$ and $t$ are odd.
\end{enumerate}
\end{obs}

\end{document}